\documentclass[11pt]{article}
\usepackage{amsmath}
\usepackage{color}
\usepackage[colorlinks=true]{hyperref}
\usepackage{url}

\def\cX{{\cal X}}
\def\cN{{\cal N}}

\newcommand{\half}{ \mbox{\small$\frac{1}{2}$}}

\usepackage{amsfonts}
\usepackage{epsfig}

\def\cZ{{\cal Z}}

\usepackage{amssymb}

\usepackage{graphicx}

\def\bR{{\mathbf{R}}}

\def\cP{{\cal P}}

\def\cX{{\cal X}}

\def\cN{{\cal N}}

\definecolor{MyDarkBlue}{rgb}{0,0.08,0.45}
\definecolor{MyViolet}{rgb}{0.45,0.08,0.95}
\definecolor{MyBrown}{rgb}{0.45,0.08,0}

\def\bE{{\mathbf{E}}}

\newtheorem{lemma}{Lemma}[section]

\newtheorem{proposition}{Proposition}[section]

\def\qed{\ \hfill$\square$\par\smallskip}

\def\e{{\rm e}}

\newcommand{\aic}[2]{{\color{MyDarkBlue}~#2}}

\newcommand{\be}{\begin{eqnarray}}
\newcommand{\ee}[1]{\label{#1}\end{eqnarray}}
\newcommand{\nn}{\nonumber \\}
\newcommand{\ese}{\end{eqnarray*}}
\newcommand{\bse}{\begin{eqnarray*}}
\newcommand{\rf}[1]{~(\ref{#1})}

\definecolor{MyDarkBlue}{rgb}{0,0.08,0.45}
\definecolor{MyViolet}{rgb}{0.45,0.08,0.95}
\definecolor{MyBrown}{rgb}{0.45,0.08,0}

\def\sML{{\textrm{\tiny ML}}}
\def\VI{{\textrm{VI}}}
\def\Proj{{\mathrm{Proj}}}
\def\argmin{\mathop\mathrm{argmin}}
\title{
Signal recovery by Stochastic Optimization}
\author{
Anatoli Juditsky
\thanks{LJK, Universit\'e Grenoble Alpes, 700 Avenue Centrale 38401 Domaine Universitaire
de Saint-Martin-d'H\`{e}res, France,
{\tt anatoli.juditsky@univ-grenoble-alpes.fr}}
\and Arkadi Nemirovski
\thanks{ISyE, Georgia Institute
 of Technology, Atlanta, Georgia
30332, USA, {\tt nemirovs@isye.gatech.edu}\newline
The first author was supported by the PGMO grant 2016-2032H. Research of
the second author was supported by NSF grant   CCF-1523768.}}
\date{}
\begin{document}
\maketitle
\begin{abstract}
We discuss an approach to signal recovery in Generalized Linear Models (GLM) in which the signal estimation problem is reduced to the problem of solving a stochastic monotone variational inequality (VI). The solution to the stochastic VI can be found in a computationally efficient way, and in the case when the VI is strongly monotone we derive finite-time upper bounds on the expected $\|\cdot\|_2^2$ error converging to 0 at the rate $O(1/K)$ as the number $K$ of observations grows.
Our structural assumptions are essentially weaker than those necessary to ensure convexity of the optimization problem resulting from Maximum Likelihood estimation. In hindsight, the approach we promote can be traced back directly to the ideas behind the Rosenblatt's perceptron algorithm.
\end{abstract}
\section{Introduction}\label{sec:intro}
Statistical estimation problems constitute one of principal application domains of Stochastic Optimization. A typical setting is as follows (cf., e.g., \cite{devroye2013probabilistic} and references therein): we are given {i.i.d.} observations $\omega^K=(\omega_1,...,\omega_K)$, $\omega_k=(\eta_k,y_k)$, where $\eta_k\in\bR^{n\times m}$, $y_k\in\bR^m$ are, respectively, realizations of regressors (independent variables) and responses (labels). {We assume that the observations can be described by a {\em Generalized Linear Model} (GLM) \cite{nelder1972generalized,mccullagh1989generalized}, that is, the conditional, $\eta$ given, expectation of $y$ is $f(\eta^Tx)$, where $f(\cdot):\bR^m\to\bR^m$ is a known {\em link function},  and $x\in\bR^m$ is the unknown ``signal'' --- vector of model's parameters. Our goal is to ``fit the model,'' that is, to recover $x$ from observations $\omega^K$.}  The standard approach to fitting the model is to choose a {\em loss function} $\ell(y,\theta):\,\bR^{m}\times \bR^m\to\bR$ {and to recover $x$ as} an optimal solution to the optimization problem
\be
{\min_{u\in \cX}\bE_{\omega\sim P_x}\{\ell(y,f(\eta^Tu))\},}
\ee{eq:stocho0}
{where} $P_x$ {is} the distribution of the observation $\omega=(\eta,y)$ associated with ``true signal $x$'', and $\cX$ is an {\em a priori} known signal set. In other words, in the just presented framework, the statistical estimation problem reduces to the stochastic optimization problem \rf{eq:stocho0}, which is to be solved {approximately via} available observation $\omega^K$. {This can be done} either ``in batch,'' {minimizing in $u\in \cX$ the} {\em Sample Average Approximation} (SAA)
{
 \begin{equation}\label{ait2}
 {1\over K}\sum_{k=1}^K\ell(y_k,f(\eta_k^Tu))
 \end{equation}}
 of the expectation in \rf{eq:stocho0} (see, e.g. \cite{ShapiroSAA}), or applying iterative stochastic optimization algorithms of {\em Stochastic Approximation} (SA) type \cite{robbins1951stochastic,wolfowitz1952stochastic}.\par
Assuming that the conditional, given $\eta$,  distribution $P^x_{|\eta}$ of $y$ induced by $P_x$ belongs to a known parameteric family $\cP=\{P^\theta:\,
\theta\in \Theta\subset\bR^m\}$, specifically, $P^x_{|\eta}=P^{f(\eta^Tx)}$, the standard choice of the loss function is given by Maximum Likelihood: assuming that distributions $P^\theta$ have densities $p_{\theta}$ w.r.t. a reference measure $\Pi$, one uses
$$
\ell(y,\theta)=-\ln(p_\theta(y)).
$$
For example, in the classical {\sl logistic regression} $m=1$, $f(s)=(1+\e^{-s})^{-1}$, $\Theta=(0,1)$, and $P^\theta$, $\theta \in\Theta$, is Bernoulli distribution, that is, label  $y$ takes value 1 with probability $(1+\exp\{-\eta^Tx\})^{-1}$ and value 0 with the complementary probability, resulting in
\be
\ell(y,f(\eta^Tu))=\ln(1+\exp\{\eta^Tu\})-y\eta^Tu.
\ee{logiteqm1}
In this case, problem (\ref{eq:stocho0}) becomes the optimization  problem
\begin{equation}\label{ait1}
\min_{u\in \cX} \bE_{(\eta,y)\sim P_x}\left\{\ln(1+\exp\{\eta^Tu\})-y\eta^Tu\right\},
\end{equation}
and its SAA  becomes
\begin{equation}\label{ait3}
\min_{u\in \cX} {1\over K}\sum_{k=1}^K\left[\ln(1+\exp\{\eta_k^Tu\})-y_k\eta_k^Tu\right];
\end{equation}
{the optimal solution} $\widehat{x}_\sML(\omega^K)$ {to the latter problem is the} {\em Maximum Likelihood (ML) estimate of $x$}.
Assuming the signal set $\cX$ to be convex, both these problems turn out to be convex, implying the possibility to solve the SAA to global optimality in a computationally efficient fasion, same as utilizing nice convergence properties of SA.
\par
\aic{}{More generally, when distributions of observations form {\em a conditional exponential family} \cite{barndorff1978information,feigin1981conditional},
negative log-likelihood has the form
\[
\{\ell(y,\eta^Tu)=F(\eta^Tu)-y\eta^Tu,
\]
with convex {\em cumulant function} $F$,
and corresponding risk minimization problem
\rf{eq:stocho0} reads
\be
\min_{u\in \cX} \bE_{(\eta,y)\sim P_x}
\left\{F(\eta^Tu)-y\eta^Tu\right\}.
\ee{eq:expf}
In this case, same as in the case of logistic regression, SAA or SA can be applied to compute Maximal Likelihood estimates of model parameter.
}
Note, however, that \aic{}{exponential family assumption is quite restrictive. On the other hand, beyond exponential families,} the convexity of the optimization problem resulting from Maximum Likelihood selection of $\ell(\cdot)$ appears to be an exception  rather than a rule. For example, consider the ``nonlinear Least Squares'' setting in which the label $y$ is obtained from $f(\eta^T x)$ by adding independent of the regressor zero mean Gaussian noise:
$$
y=f(\eta^Tx)+\xi,\,\,\xi\sim\cN(0,\sigma^2I_m).
$$
In this case problem (\ref{eq:stocho0}) and its SAA approximation for the ML selection of $\ell(\cdot)$ become
\be
&\min_{u\in \cX}\bE_{\eta\sim Q}\left\{\|f(\eta^Tx)-f(\eta^Tu)\|_2^2\right\},\nn
&\min_{u\in \cX}\left\{{1\over K}\sum_{k=1}^k\|y_k-f(\eta_k^Tu)\|_2^2\right\},
\ee{MLnot}
where $Q$ is the distribution of regressors (which we assume to be independent of the signal). When $f$ is nonlinear, both these problems usually are nonconvex and could be difficult to process numerically. Similarly, in the ``non-exponential logistic regression,''
where the ``exponential sigmoid'' $f(s)=(1+\exp\{-s\})^{-1}$ is replaced with a general nondecreasing link function $f(s):\bR\to(0,1)$ (e.g., {\em probit} or {\em complementary $\log-\log$} link) 
the ML selection of the loss function typically makes (\ref{eq:stocho0}) and its SAA approximation nonconvex.
\par
The goal of what follows is to propose an alternative to  model fitting via (\ref{eq:stocho0}) with ML-based selection of the loss function approach to estimating the signal underlying observations in a GLM.
In hindsight, the approach we put forward in this paper can be traced back to the ideas behind the Rosenblatt's perceptron iterative algorithm \cite{rosenblatt1958perceptron,block1962perceptron} and its batch version \cite{helmbold1995weak}.  The structural assumptions to be imposed on the model are essentially weaker than those resulting in convex {ML-based} problems (\ref{eq:stocho0}) and their SAA approximations.\footnote{For instance, in the ``nonlinear least squares'' with $m=1$, same as in ``non-exponential logistic regression,'' all we need from $f$ to be continuously differentiable, with positive derivative, and from the signal set $\cX$ to be convex.} Under these assumptions, instead of using the classical loss function approach \cite{aiserman1964theoretical,devyaterikov1967iterative,aizerman1970method,bousquet2004introduction,sridharan2009fast},
we reduce the estimation problem to another problem with convex structure --- a {\sl strongly monotone variational inequality} (VI) {\sl represented by a stochastic oracle.} {This VI may or may not be equivalent to a convex minimization problem. The first option definitely takes place when $m=1$,
\aic{}{when the VI is equivalent to analogous to \rf{eq:expf} convex optimization problem;} but even in this case the resulting problem typically is different from the ML version of (\ref{eq:stocho0}).} The solution to the {VI} can be found in a computationally efficient way and turns out to be a ``good'' estimate of the signal underlying observations, for which we derive finite-time upper bounds on the expected $\|\cdot\|_2^2$ error, converging to 0 at the rate $O(1/K)$ as $K\to\infty$.\footnote{We were unable to locate a reference to the proposed approach in the statistical literature, though it would be the most surprising if simple derivations which follow were not known.}

\section{Problem statement}
Throughout the paper we consider the GLM model as posed in Introduction:
\begin{quote}
Our observation depends on unknown signal $x$ known to belong to a given convex compact set $\cX\subset\bR^n$ and is
\begin{equation}\label{logiteq1}
\omega^K=\{\omega_k=(\eta_k,y_k),1\leq k\leq K\}
\end{equation}
with  $\omega_k$, $1\leq k\leq K$ which are i.i.d. realizations of a random pair $(\eta,y)$ with the distribution $P_x$ such that
\begin{itemize}
\item the {\sl regressor} $\eta$ is a random $n\times m$ matrix with some {\sl independent of $x$} probability distribution $Q$;
\item the {\sl label} $y$ is $m$-dimensional random vector such that the conditional, given $\eta$, distribution
of $y$ induced by $P_x$ has the expectation $f(\eta^Tx)$:
\begin{equation}\label{logiteq2}
\bE^x_{|\eta}\{y\}=f(\eta^Tx),
\end{equation}
where
 $\bE^x_{|\eta}$ is the conditional, $\eta$ given, distribution of $y$ stemming from the distribution $P_x$ of $\omega=(\eta,y)$, and
$f(\cdot):\bR^m\to\bR^m$ is a given mapping.
\end{itemize}
\end{quote}
We are about to formulate assumptions on the parameters of a generalized linear model (namely, on $f(\cdot)$, and the distributions $P_x$, $x\in\cX$, of the pair $(\eta,y)$) required by the approach we are about to develop.
\subsection{Preliminaries: monotone vector fields} A  monotone vector field on $\bR^m$ is a single-valued everywhere defined mapping $g(\cdot):\bR^m\to\bR^m$ which  possesses the {\sl monotonicity property}
$$
[g(z)-g(z')]^T[z-z']\geq 0\,\,\forall z,z'\in\bR^m.
$$
We say that such a field is {\sl monotone with modulus $\varkappa\geq0$  on a closed convex set $Z\subset\bR^m$}, if
$$
[g(z)-g(z')]^T[z-z']\geq \varkappa\|z-z'\|_2^2\\,\forall z\,z'\in Z,
$$
and say that $g$ is {\sl strongly monotone} on $Z$ if the modulus of monotonicity of $g$ on $Z$ is positive.
It is immediately seen that for a monotone vector field which is continuously differentiable on a closed convex set $Z$ with a nonempty interior, the necessary and sufficient condition for being monotone with modulus $\varkappa$ on the set is
\begin{equation}\label{logiteqm2}
d^Tf'(z)d\geq\varkappa d^Td\,\,\forall (d\in\bR^n,z\in Z).
\end{equation}
Basic examples of monotone vector fields are:
\begin{itemize}
\item gradient fields $\nabla \phi(x)$ of continuously differentiable convex functions of $m$ variables or, more generally, the
vector fields $[\nabla_x \phi(x,y);-\nabla_y\phi(x,y)]$ stemming from continuously differentiable functions $\phi(x,y)$ which are convex in $x$ and concave in $y$;
\item ``diagonal'' vector fields $f(x)=[f_i(x_1);f_2(x_2);...;f_m(x_m)]$ with monotonically nondecreasing univariate components $f_i(\cdot)$.
If, in addition, $f_i(\cdot)$ are continuously differentiable with positive derivatives, then the associated field $f$ is strongly monotone on every compact convex subset of $\bR^m$, the monotonicity modulus depending on the subset.
\end{itemize}
Monotone vector fields on $\bR^n$ admit simple calculus which includes, in particular, the following two rules:
\begin{itemize}
\item[\textbf{I}.] [affine substitution of argument]: If $f(\cdot)$ is monotone vector field on $\bR^m$ and $A$ is an $n\times m$ matrix, the vector field
$$
g(x)=Af(A^Tx+a)
$$ is monotone on $\bR^n$; if, in addition, $f$ is monotone with modulus $\varkappa\geq0$ on a closed convex set $Z\subset\bR^m$ and $X\subset\bR^n$ is closed, convex, and such that $A^Tx+a\in Z$ whenever $x\in X$, $g$ is monotone with modulus $\sigma^2\varkappa$ on $X$,  where $\sigma$ is the minimal singular value of $A$.
\item[\textbf{II}.] [summation]: If $S$ is a Polish  space, $f(x,s):\bR^m\times S\to\bR^m$ is a Borel vector-valued function which is monotone in $x$ for
every $s\in S$ and $\mu(ds)$ is a Borel probability measure on $S$ such that the vector field
$$
F(x)=\int_Sf(x,s)\mu(ds)
$$
is well defined for all $x$, then $F(\cdot)$ is monotone. If, in addition, $X$ is a closed convex set in $\bR^m$ and $f(\cdot,s)$ is monotone on $X$ with Borel in $s$ modulus $\varkappa(s)$ for every $s\in S$, then $F$ is monotone on $X$ with modulus $\int_S\varkappa(s)\mu(ds)$.
\end{itemize}
\subsection{Assumptions} In what follows, we make the following assumptions on the ingredients of the estimation problem set in
Introduction:
\begin{itemize}
\item{\textbf{A.1}.} The vector field $f(\cdot)$ is continuous and monotone, and the vector field
$$
F(z)=\bE_{\eta\sim Q}\left\{\eta f(\eta^Tz)\right\}
$$
is well defined (and therefore is monotone along with $f$ by \textbf{I}, \textbf{II});
\item{\textbf{A.2}.} The signal set $\cX$ is a nonempty convex compact set, and the vector field $F$  is monotone with positive modulus $\varkappa$ on $\cX$;
\item{\textbf{A.3}.} For properly selected $M<\infty$ and every $x\in\cX$ it holds
\begin{equation}\label{logiteq3}
\bE_{(\eta,y)\sim P_x}\left\{\|\eta y\|_2^2\right\}\leq  M^2.
\end{equation}
\end{itemize}
A simple {\sl sufficient} condition for the validity of Assumptions \textbf{A.1-3} with properly selected $M<\infty$ and $\varkappa>0$ is as follows:
\begin{itemize}
\item The distribution $Q$ of $\eta$ has finite moments of all orders, and $\bE_{\eta\sim Q} \{\eta\eta^T\}\succ0$;
\item $f$ is continuously differentiable, and $d^Tf'(z)d>0$ for all $d\neq0$ and all $z$. Besides this, $f$ is of polynomial growth: for some constants $C\geq0$ and $p\geq0$ and all $z$ one has $\|f(z)\|_2\leq C(1+\|z\|_2^p)$.
    \end{itemize}
Verification of sufficiency is straightforward.

\section{Construction and Main result}\label{logitmainobs}
The principal observation underlying the construction we are about to present is as follows:
\begin{proposition}\label{logitobs1}
Assuming that Assumptions \textbf{A.1-3} hold, let us associate with the pair $(\eta,y)\in\bR^{n\times m}\times\bR^m$
the vector field
\begin{equation}\label{logiteq4}
G_{(\eta,y)}(z)=\eta f(\eta^Tz)-\eta y:\bR^n\to\bR^n.
\end{equation}
Then for every $x\in\cX$
we have
\begin{equation}\label{logiteq5}
\begin{array}{rclr}
 \bE_{(\eta,y)\sim P_x}\left\{G_{(\eta,y)}(z)\right\}&=&F(z)-F(x)\,\,\forall z\in\bR^n&(a)\\
 \|F(z)\|_2&\leq& M\,\forall z\in\cX&(b)\\
 \bE_{(\eta,y)\sim P_x}\left\{\|G_{(\eta,y)}(z)\|_2^2\right\}&\leq& 4M^2\,\,\forall z\in\cX&(c)\\
 \end{array}
\end{equation}
\end{proposition}
\textbf{Proof} is immediate. Indeed, let $x\in\cX$. Then
$$
\bE_{(\eta,y)\sim P_x} \{\eta y\}=\bE_{\eta\sim Q}\left\{\bE^x_{|\eta}\{\eta y\}\right\}=
\bE_{\eta}\left\{\eta f(\eta^Tx)\right\}=F(x)
$$
(we have used (\ref{logiteq2}) and the definition of $F$), whence,
\bse
\bE_{(\eta,y)\sim P_x}\left\{G_{(\eta,y)}(z)\right\}&=&\bE_{(\eta,y)\sim P_x}\left\{\eta f(\eta^T z)-\eta y\right\}
=\bE_{(\eta,y)\sim P_x}\left\{\eta f(\eta^T z)\right\}-F(x)\\
&=&\bE_{\eta\sim Q}\left\{\eta f(\eta^T z)\right\}-F(x)=F(z)-F(x),
\ese as stated in (\ref{logiteq5}.$a$). Besides this, for $x,z\in \cX$, denoting by $P^z_{|\eta}$ the conditional, $\eta$ given, distribution of $y$ induced by the distribution $P_z$ of $(\eta,y)$, and taking into account that the marginal distribution of $\eta$ induced by $P_z$ is $Q$, we have
\bse
\bE_{(\eta,y)\sim P_x}\{\|\eta f(\eta^T z)\|_2^2\}&=&\bE_{\eta\sim Q} \left\{\|\eta f(\eta^Tz)\|_2^2\right\}\\
&=&\bE_{\eta\sim Q}\left\{\|\bE_{y\sim P^z_{|\eta}}\{\eta y\}\|_2^2\right\} \hbox{\ [since $\bE_{y\sim P^z_{|\eta}}\{y\}=f(\eta^Tz)$]}\\
&\leq& \bE_{\eta\sim Q}\left\{\bE_{y\sim P^z_{|\eta}}\left\{\|\eta y\|_2^2\right\}\right\}\hbox{\ [by Jensen's inequality]}\\
&=&\bE_{(\vec{}\eta,y)\sim P_z}\left\{\|\eta y\|_2^2\right\}
\leq M^2\hbox{\ [by \textbf{A.3} due to $z\in\cX$]}.
\ese
This combines with the relation $\bE_{(\eta,y)\sim P_x} \{\|\eta y\|_2^2\}\leq M^2$ given by \textbf{A.3} due to $x\in\cX$ to imply
(\ref{logiteq5}.$b$) and (\ref{logiteq5}.$c$). \qed
\subsection{Main result}
Recall that our goal is to recover the signal $x\in\cX$ underlying observations (\ref{logiteq1}). Under assumptions \textbf{A.1-3}, $x$ is a root of the monotone vector field
\begin{equation}\label{logiteq44}
G(z)=F(z)-F(x),\,\,F(z)=\bE_{\eta\sim Q}\left\{\eta f(\eta^Tz)\right\};
\end{equation} we know that this root belongs to $\cX$, and is unique because $G(\cdot)$ is strongly monotone on $\cX$ along with $F(\cdot)$. Now, finding a root, known to belong to a given convex compact set $\cX$,  of a strongly monotone on this set vector field $G$ is known to be a computationally tractable problem, provided we have access to an ``oracle'' which, given on input a point $z\in \cX$, returns the value $G(z)$  of the field at the point.  The latter is not exactly the case in the situation we are interested in: the field $G$ is the expectation of a random field:
$$
G(z)=\bE_{(\eta,y)\sim P_x}\left\{\eta f(\eta^Tz)-\eta y\right\},
$$
and we do not know a priori what is the distribution over which the expectation is taken. However, we can sample from this distribution -- the samples are exactly the observations (\ref{logiteq1}), and we can use these samples to approximate somehow $G$ and use this approximation to approximate the signal $x$. Two standard implementations of this idea are {\sl Sample Average Approximation} (SAA) and {\sl Stochastic Approximation} (SA). We are about to consider these two techniques as applied to the situation we are in.
\index{Sample Average Approximation|(}

\subsubsection{Estimation by Sample Average Approximation}
The idea underlying SAA is quite transparent: given observations (\ref{logiteq1}), let us approximate the field of interest $G$ with its  empirical counterpart
$$
G_{\omega^K}(z)={1\over K}\sum_{k=1}^K \left[\eta_kf(\eta_k^Tz)-\eta_ky_k\right].
$$ By the Law of Large Numbers, as $K\to\infty$, the empirical field $G_{\omega^K}$ converges to the field of interest $G$, so that under mild regularity assumptions, when $K$ is large, $G_{\omega^K}$, with overwhelming  probability, will be uniformly on $\cX$ close to $G$.
Due to strong monotonicity of $G$, this would imply that a set of ``near-zeros'' of $G_{\omega^K}$ on $\cX$ will be close to the zero $x$ of $G$, which is nothing but the signal we want to recover. The only question is how we can consistently define a
``near-zero'' of $G_{\omega^K}$ on $\cX$.\footnote{Note that we in general cannot define a ``near-zero'' of $G_{\omega^K}$ on $\cX$ as a root of $G_{\omega^K}$ on this set -- while $G$ does have a root belonging to $\cX$, nobody told us that the same holds true for $G_{\omega^K}$.} A convenient in our context notion of a ``near-zero'' is provided by the concept of a {\sl weak solution} to a variational inequality (VI) with monotone operator, defined as follows (we restrict the general definition to the situation of interest):
\begin{quote}
Let $\cX\subset\bR^n$ be a nonempty convex compact set, and $H(z): \cX\to\bR^n$ be a monotone (i.e., $[H(z)-H(z')]^T[z-z']\geq0$ for all $z,z'\in \cX$) vector field. A vector $z_*\in\cX$ is called a {\sl weak solution} to the variational inequality (VI) associated with $H,\cX$ when
$$
H(z)^T(z-z_*)\geq0\;\forall z\in\cX.
$$
\end{quote}
Let $\cX\subset\bR^n$ be a nonempty convex compact set and $H$ be monotone on $\cX$. It is well known that
\begin{itemize}
\item The VI associated with $H,\cX$ (let us denote it $\VI(H,\cX)$)  always has a weak solution. It is clear that if $\bar{z}\in\cX$ is a root of $H$, then $\bar{z}$ is a weak solution to $\VI(H,\cX)$.\footnote{Indeed, when $\bar{z}\in\cX$ and $H(\bar{z})=0$, monotonicity of $H$ implies  that $H(z)^T[z-\bar{z}]=[H(z)-H(\bar{z})]^T[z-\bar{z}]\geq0$ for all $z\in\cX$, that is, $\bar{z}$ is a weak solution to the VI.}
\item  When $H$ is continuous on $\cX$, every weak solution $\bar{z}$ to $\VI(H,\cX)$ is also a {\sl strong solution}, meaning that
\begin{equation}\label{logiteq60}
H^T(\bar{z})(z-\bar{z})\geq0\,\,\forall z\in\cX.
\end{equation}
Indeed, (\ref{logiteq60}) clearly holds true when $z=\bar{z}$. Assuming $z\neq\bar{z}$ and setting $z_t=\bar{z}+t(z-\bar{z})$, $0<t\leq 1$, we have $H^T(z_t)(z_t-\bar{z})\geq0$ (since $\bar{z}$ is a weak solution), whence $H^T(z_t)(z-\bar{z})\geq0$ (since $z-\bar{z}$ is a positive multiple of $z_t-\bar{z}$). Passing to limit as $t\to+0$ and invoking the continuity of $H$, we get $H^T(\bar{z})(z-\bar{z})\geq0$, as claimed.
\item When $H$ is the gradient field of a continuously differentiable convex function on $\cX$ (such a field indeed is monotone),  weak (or,  which in the case of continuous $H$ is the same, strong) solutions to $\VI(H,\cX)$ are exactly the minimizers of the function on $\cX$.\par
Note also that a strong solution to $\VI(H,\cX)$ with monotone $H$ always is a weak one: if $\bar{z}\in\cX$ satisfies $H^T(\bar{z})(z-\bar{z})\geq0$ for all $z\in\cX$, then $H(z)^T(z-\bar{z})\geq0$ for all $z\in \cX$, since by monotonicity $H(z)^T(z-\bar{z})\geq H^T(\bar{z})(z-\bar{z})$.
\end{itemize}
In the sequel, we heavily exploit the following simple and well known fact:
{\begin{lemma}
\label{logitlem} Let $\cX$ be a convex compact set, and $H$ be a monotone vector field on $\cX$ with monotonicity modulus $\varkappa>0$, i.e.
\[
\forall z,z'\in X\;[H(z)-H(z')]^T[z-z']\geq \varkappa\|z-z'\|_2^2.
\]
Further, let  $\bar{z}$ be a weak solution to $\VI(H,\cX)$. Then
the  weak solution to $\VI(H,\cX)$ is unique. Besides this,
\begin{equation}\label{logiteq8}
H^T(z)[z-\bar{z}]\geq \varkappa \|z-\bar{z}\|_2^2.
\end{equation}
\end{lemma}}\noindent
\textbf{Proof:} Under the premise of the lemma, let $z\in \cX$ and let $\bar{z}$ be a weak solution to $\VI(H,\cX)$ (recall that it does exist). Setting $z_t=\bar{z}+t(z-\bar{z})$, for $t\in(0,1)$ we have
$$
H^T(z)[z-z_t]\geq H^T(z_t)[z-z_t]+\varkappa \|z-z_t\|^2\geq \varkappa \|z-z_t\|^2,
$$
where the first $\geq$ is due to strong monotonicity of $H$, and the second $\geq$ is due to the fact that $H^T(z_t)[z-z_t]$ is proportional, with positive coefficient, to $H^T(z_t)[z_t-\bar{z}]$, and the latter quantity is nonnegative since $\bar{z}$ is a weak solution to the VI in question. We end up with $H^T(z)(z-z_t)\geq \varkappa\|z-z_t\|_2^2$; passing to limit as $t\to+0$, we arrive at (\ref{logiteq8}). To prove uniqueness of a weak solution, assume that aside of the weak solution $\bar{z}$ there exists a weak solution $\widetilde{z}$ distinct form $\bar{z}$, and let us set $z'={1\over 2}[\bar{z}+\widetilde{z}]$. Since both $\bar{z}$ and $\widetilde{z}$ are weak solutions, both the quantities $H^T(z')[z'-\bar{z}]$
and $H^T(z')[z'-\widetilde{z}]$ should be nonnegative, and because the sum of these quantities is 0, both of them are zero. Thus, when applying (\ref{logiteq8}) to $z=z'$, we get $z'=\bar{z}$, whence $\widetilde{z}=\bar{z}$ as well. \qed

Now, let us return to the estimation problem under consideration. Assume that Assumptions \textbf{A.1-3} hold, so vector fields  $G_{(\eta_k,y_k)}(z)$ defined in (\ref{logiteq4}), and therefore vector field $G_{\omega^K}(z)$ are continuous and monotone. When using the SAA, we compute a weak solution $\widehat{x}(\omega^K)$ to $\VI(G_{\omega^K},\cX)$ and treat it as the SAA estimate of signal $x$ underlying observations (\ref{logiteq1}). Since the vector field $G_{\omega^K}(\cdot)$ is monotone with  efficiently computable values, provided that so is $f$, computing (a high accuracy approximation to) a weak solution to $\VI(G_{\omega^K},\cX)$ is a computationally tractable problem (see, e.g., \cite{NOR}).
Moreover, utilizing the techniques from \cite{bousquet2002stability,rakhlin2005stability,sridharan2009fast,shalev2009stochastic}, under mild additional to \textbf{A.1-3} regularity assumptions one can get non-asymptotical upper bound on, say, the expected $\|\cdot\|_2^2$-error of the SAA estimate as a function of the sample size $K$ and find out the rate at which this bound converges to 0 as $K\to\infty$; this analysis, however, goes beyond our scope.
\par
Let us {look at} the SAA estimate in the logistic regression model. In this case we have $f(u)=(1+\e^{-u})^{-1}$, and
\bse
G_{(\eta_k,y_k)}(z)&=&\left[{\exp\{\eta_k^Tz\}\over 1+\exp\{\eta_k^Tz\}}-y_k\right]\eta_k,\\
G_{\omega^K}(z)&=&{1\over K}\sum_{k=1}^K\left[{\exp\{\eta_k^Tz\}\over 1+\exp\{\eta_k^Tz\}}-y_k\right]\eta_k\\
&=&{1\over K}\nabla_z\left[\sum_k\left(\ln\left(1+\exp\{\eta_k^Tz\}\right)-y_k\eta_k^Tz\right)\right].
\ese
In other words, $G_{\omega^K}(z)$ is  the gradient field of the minus empirical log-likelihood $\ell(z,\omega^K)$, see
(\ref{ait3}). As a result, in the case in question  weak solutions to $\VI(G_{\omega^K},\cX)$ are exactly the optimal solutions to (\ref{ait3}), that is, {\sl for the logistic regression the SAA estimate is nothing but the Maximum Likelihood estimate $\widehat{x}_{\sML}(\omega^K)$.}
\footnote{This phenomenon is specific for the logistic regression model. The fact that the SAA and the ML estimates in this case are the same is due to the fact
that the logistic sigmoid $f(s)=\exp\{s\}/(1+\exp\{s\})$ ``happens'' to satisfy the identity $f'(s)=f(s)(1-f(s))$. When replacing the exponential sigmoid with $f(s)=\phi(s)/(1+\phi(s))$ with differentiable monotonically nondecreasing positive $\phi(\cdot)$, the SAA estimate becomes the weak solution to $\VI(\Phi,\cX)$ with
$$
\Phi(z)=\sum_k \left[{\phi(\eta_k^Tz)\over 1+\phi(\eta_k^Tz)}-y_k\right]\eta_k.
$$
On the other hand, the gradient field of the {\sl minus} log-likelihood
$-{1\over K}\sum_k\left[y_k\ln(f(\eta_k^Tz))+(1-y_k)\ln(1-f(\eta_k^Tz))\right]$) which we should minimize when computing the ML estimate
is
$$\Psi(z)=\sum_k{\phi'(\eta_k^Tz)\over\phi(\eta_k^Tz)}\left[{\phi(\eta_k^Tz)\over 1+\phi(\eta_k^Tz)}-y_k\right]\eta_k.
$$
When $k>1$ and $\phi$ is not an exponent, $\Phi$ and $\Psi$ are ``essentially different,'' so that the SAA estimate typically will differ from the ML one.}
 On the other hand, in the ``nonlinear least squares'' example described in the introduction with (for the sake of simplicity, scalar) monotone $f(\cdot)$
the vector field $G_{\omega^K}(\cdot)$ is given by
$$
G_{\omega^K}(z)={1\over K}\sum_{k=1}^K \left[f(\eta_k^Tz)-y_k\right]\eta_k
$$
which is ``essentially different'' (provided that $f$ is nonlinear) from the gradient field
$$
\Psi(z)={2\over K}\sum_{k=1}^Kf'(\eta_k^Tz)\left[f(\eta_k^Tz)-y_k\right]\eta_k
$$
of the negative log-likelihood appearing in (\ref{MLnot}). As a result, in this case the ML estimate (\ref{MLnot}) is, in general, different from the SAA estimate (and, in contrast  to the ML, the SAA estimate is easy to compute).
\index{Sample Average Approximation|)}
\index{Stochastic Approximation|(}
\subsubsection{Stochastic Approximation estimate}
The {\sl Stochastic Approximation} (SA) estimate stems from a simple algorithm -- {\sl Subgradient Descent} -- for solving variational inequality $\VI(G,\cX)$. Were the values of the vector field $G(\cdot)$ available, one could approximate a root $x\in\cX$ of this VI using the recurrence
$$
z_{k}=\Proj_{\cX}[z_{k-1}-\gamma_k G(z_{k-1})],\,k=1,2,...,K,
$$
where
\begin{itemize}
\item $\Proj_{\cX}[z]$ is the metric projection of $\bR^n$ onto $\cX$:
$$
\Proj_{\cX}[z]=\argmin_{u\in\cX} \|z-u\|_2;
$$
\item $\gamma_k>0$ are given stepsizes;
\item the initial point $z_0$ is an arbitrary point of $\cX$.
\end{itemize}
It is well known that under Assumptions \textbf{A.1-3} this recurrence with properly selected stepsizes and started at a point from $\cX$  allows to approximate the root of $G$ (in fact, the unique weak solution to $\VI(G,\cX)$) to a whatever high accuracy, provided $K$ is large enough. However, we are in the situation when the actual values of $G$ are not available; the standard way to cope with this difficulty is to replace in the above recurrence the ``unobservable'' values $G(z_{k-1})$ of $G$ with their unbiased random estimates $G_{(\eta_k,y_k)}(z_{k-1})$. This modification gives rise to {\sl Stochastic Approximation} (coming back to \cite{KieferSA}) -- the recurrence
\begin{equation}\label{logiteqSA}
z_{k}=\Proj_{\cX}[z_{k-1}-\gamma_k G_{(\eta_k,y_k)}(z_{k-1})],\,1\leq k\leq K,
\end{equation}
where $z_0$ is a once for ever chosen point from $\cX$, and $\gamma_k>0$ are deterministic.
\paragraph{Convergence analysis.} The following result is perfectly well known; to make the paper self-contained, we present its (completely standard) proof in Appendix.

{\begin{proposition}\label{logitprop} Under Assumptions  \textbf{A.1-3} and with the stepsizes
\begin{equation}\label{logiteq50}
\gamma_k=[\varkappa (k+1)]^{-1},\,k=1,2,...
\end{equation}
for every signal $x\in\cX$ the sequence of estimates $\widehat{x}_k(\omega^{k})=z_{k}$ 
given by the SA recurrence {\rm (\ref{logiteqSA})} and $\omega_k=(\eta_k,y_k)$ defined in {\rm (\ref{logiteq1})} for every $k$ obeys the error bound
\begin{equation}\label{logiteq51}
\bE_{\omega^k\sim P_x^k}\left\{\|\widehat{x}_k(\omega^k)-x\|_2^2\right\}\leq {4M^2\over\varkappa^2(k+1)},\,k=0,1,...
\end{equation}
$P_x$ being the distribution of $(\eta,y)$ stemming from signal $x$.
\end{proposition}}\noindent
\index{Stochastic Approximation|)}
\subsection{Numerical illustration}
\begin{figure}[h]
\begin{center}
\begin{tabular}{cc}
\includegraphics[width=0.45\textwidth]{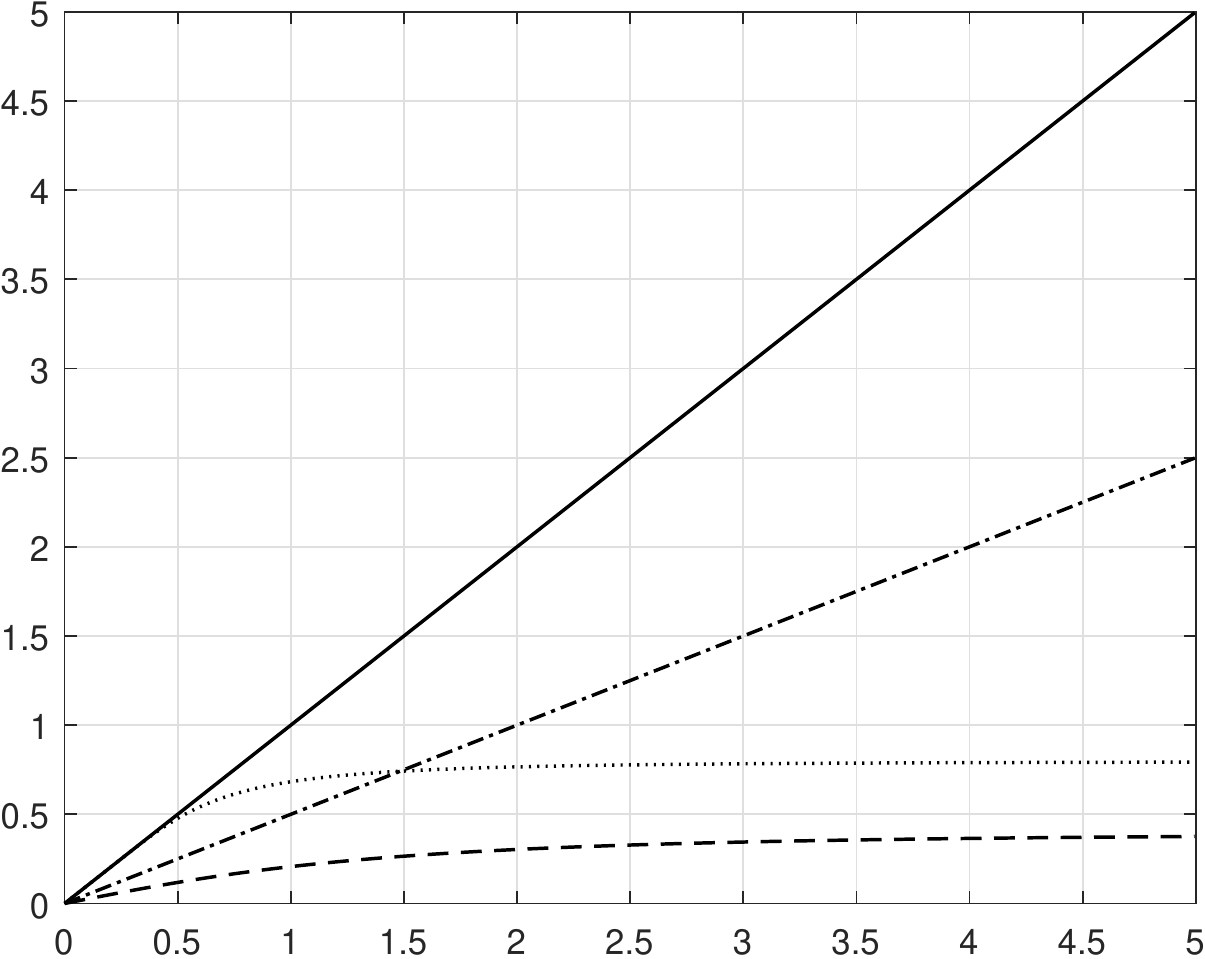}&\includegraphics[width=0.45\textwidth]{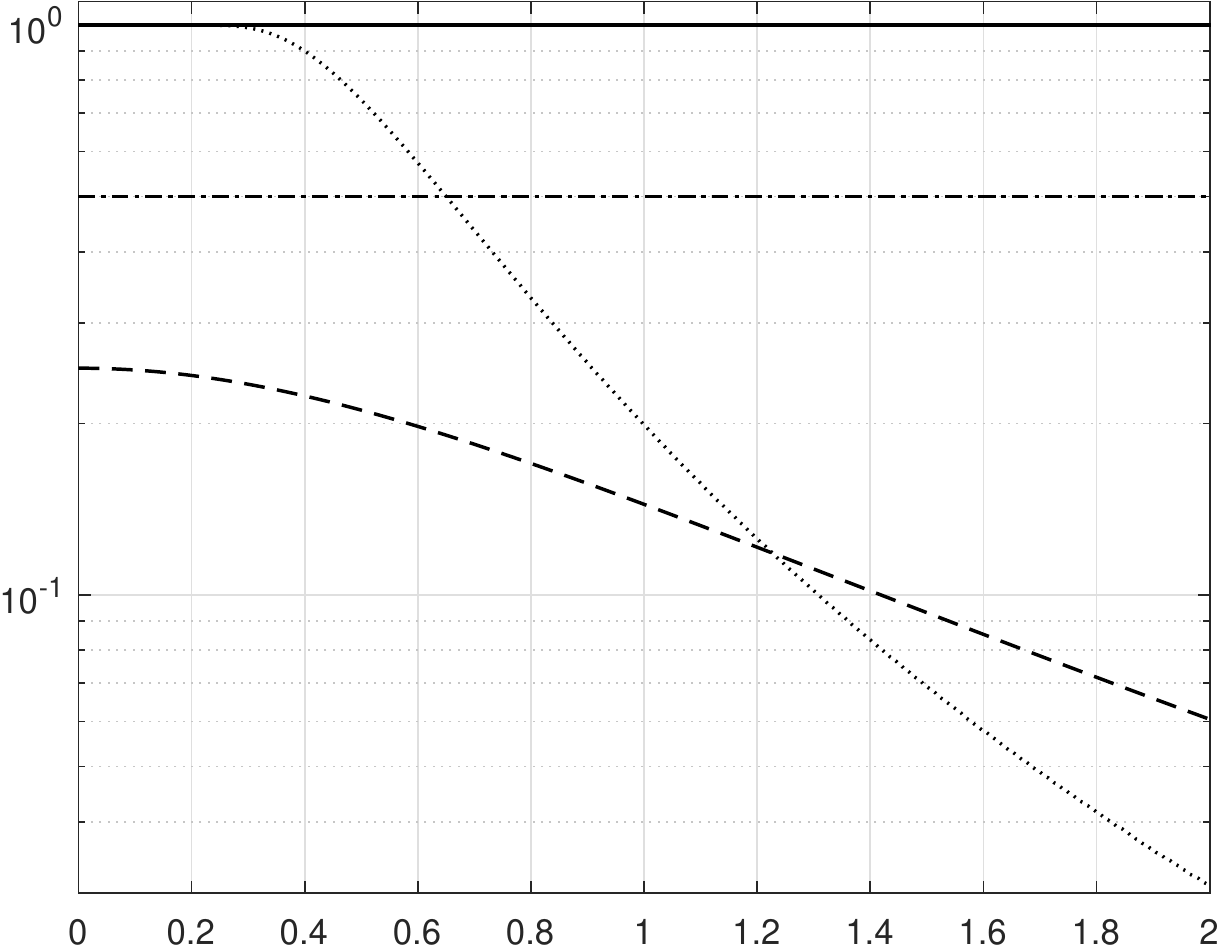}\\
\end{tabular}
\end{center}
\caption{\label{fig:funplot} Left: functions $h$; right: moduli of strong monotonicity of the operators $F(\cdot)$ in $\{z:\|z\|_2\leq R\}$ as functions of $R$.
 Dashed lines --  case A (logistic sigmoid), solid lines -- case B (linear regression), dash-dotted lines -- case C (hinge function), dotted line -- case D (``ramp'' sigmoid).}
 \end{figure}

To illustrate the above developments, we present here results of some numerical experiments. Our deliberately simplistic setup is as follows:
\begin{itemize}
\item $\cX=\{x\in\bR^n:\|x\|_2\leq1\}$;
\item the distribution $Q$ of $\eta$ is $\cN(0,I_n)$;
\item $f$ is the monotone vector field on $\bR$ given by one of the following four options:
\begin{itemize}
\item[A.] $f(s)=\exp\{s\}/(1+\exp\{s\})$;
\item[B.] $f(s)=s$;
\item[C.] $f(s)=\max[s,0]$;
\item[D.] $f(s)=\min[1,\max[s,0]]$.
\end{itemize}
\item conditional, given $\eta$, distribution of $y$ induced by $P_x$ is
\begin{itemize}
\item Bernoulli distribution with probability $f(\eta^Tx)$ of outcome 1 in the case of A (i.e., A corresponds to the logistic model),
\item Gaussian distribution $\cN(f(\eta^Tx),I_n)$ in cases B -- D.
\end{itemize}
\end{itemize}
Note that in the considered example one can easily compute the field $F(z)$. Indeed, we have $\forall z\in \bR^n$:
\[
\eta={zz^T\over \|z\|_2^2}\eta+\underbrace{\left(I-{zz^T\over \|z\|_2^2}\right)\eta}_{\eta_\perp},
\]
and due to the independence of $ \eta^Tz$ and $\eta_\perp$,
\[
F(z)=\bE_{\eta\sim \cN(0, I)}\{\eta f(\eta^Tz)\}=\bE_{\eta\sim \cN(0, I)}\left\{{zz^T\eta\over \|z\|_2^2} f(\eta^Tz)\right\}={z\over \|z\|_2}\bE_{\zeta\sim \cN(0, 1)}\{\zeta f(\|z\|_2\zeta)\},
\]
and $F(z)$ is proportional to $z/\|z\|_2$ with proportionality coefficient
$$h(\|z\|_2)=\bE_{\zeta\sim \cN(0, 1)}\{\zeta f(\|z\|_2\zeta)\}.$$
In Figure \ref{fig:funplot} we present the plots of the function $h(t)$ for the situations A -- D, same as the dependencies of the moduli of strong convexity
of the corresponding mappings $F$ in a centered at the origin $\|\cdot\|_2$-ball of radius $R$ on $R$.
The dimension $n$ in all experiments was set to 100, and the number of observations $K$ was $400$, 1e3, 4e3, 1e4, and 4e4. For each combination of parameters we ran 10 simulations for signals $x$ underlying observations (\ref{logiteq1}) drawn randomly from the uniform distribution on the unit sphere (the boundary of $\cX$).
\par
In each experiment,
we computed the SAA and the SA estimates (note that in the cases A and B  the SAA estimate is the Maximum Likelihood estimate as well).
The SA  stepsizes $\gamma_k$ were selected according to (\ref{logiteq50}) with ``empirically selected'' $\varkappa$. \footnote{We could get (lower bounds on) the modules of strong monotonicity of the vectors fields $F(\cdot)$ we are interested in analytically, but this would be boring and conservative.} Namely, given observations $\omega_k=(\eta_k,y_k)$, $k\leq K$, see (\ref{logiteq1}), we used them to build the SA estimate in two stages:
\\
--- at {\sl tuning stage}, we generate a random ``training signal'' $x'\in\cX$ and then generate labels $y_k^\prime$ as if $x'$ were the actual signal. For instance, in the case of A, $y_k^\prime$ is assigned value 1 with probability $f(\eta_k^Tx')$ and value 0 with complementary probability. After ``training signal'' and associated labels are generated, we run on the resulting artificial observations SA with different values of $\varkappa$, compute the accuracy of the resulting estimates, and select the value of $\varkappa$ resulting in the best recovery;\\
--- at {\sl execution stage}, we run SA on the actual data with stepsizes (\ref{logiteq50}) specified by $\varkappa$ found at the tuning stage.
\par
The results of some numerical experiments are presented in
{Figure \ref{fig:mono}.
\begin{figure}[h]
\begin{center}
\begin{tabular}{cc}
\includegraphics[width=0.45\textwidth]{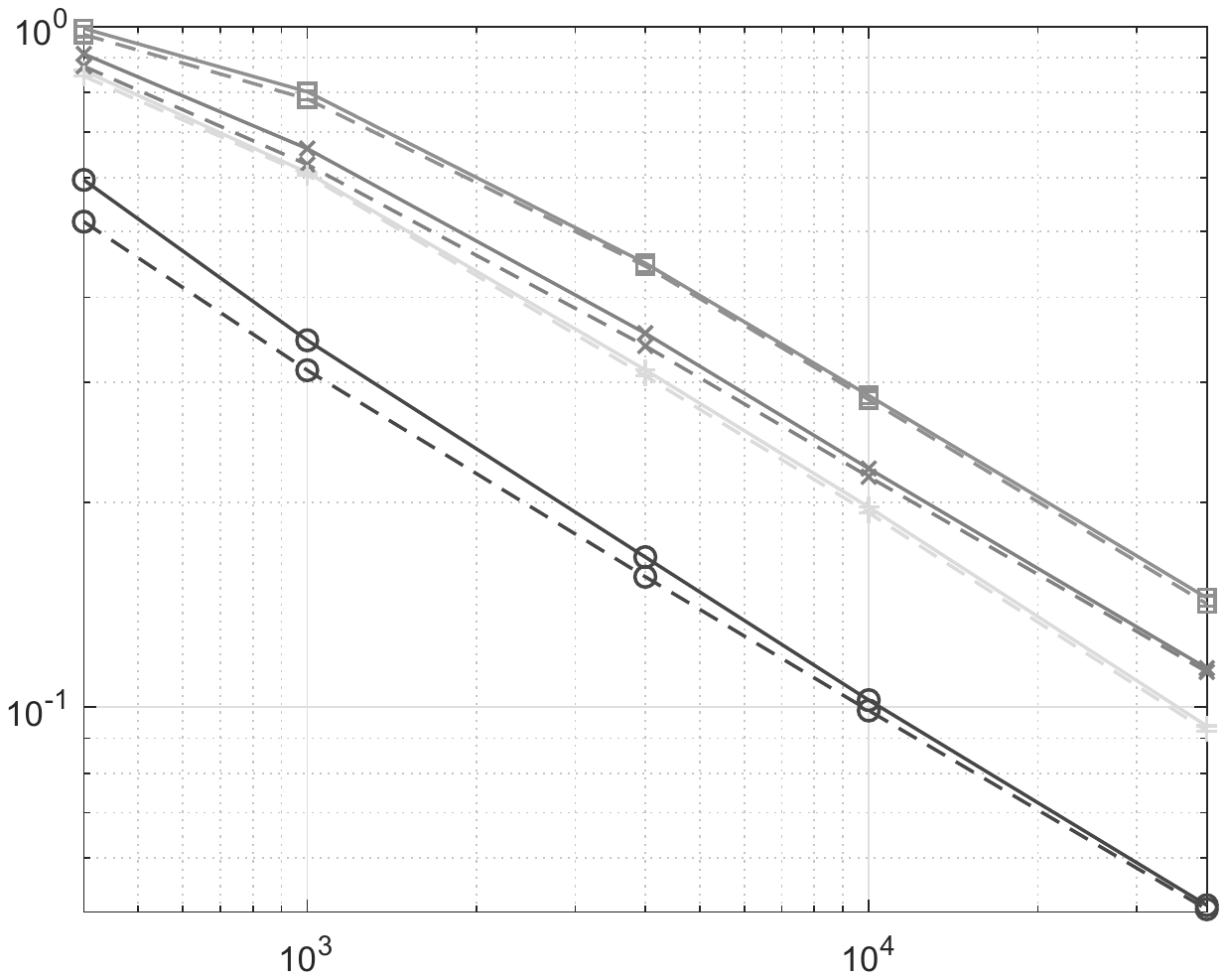}&\includegraphics[width=0.45\textwidth]{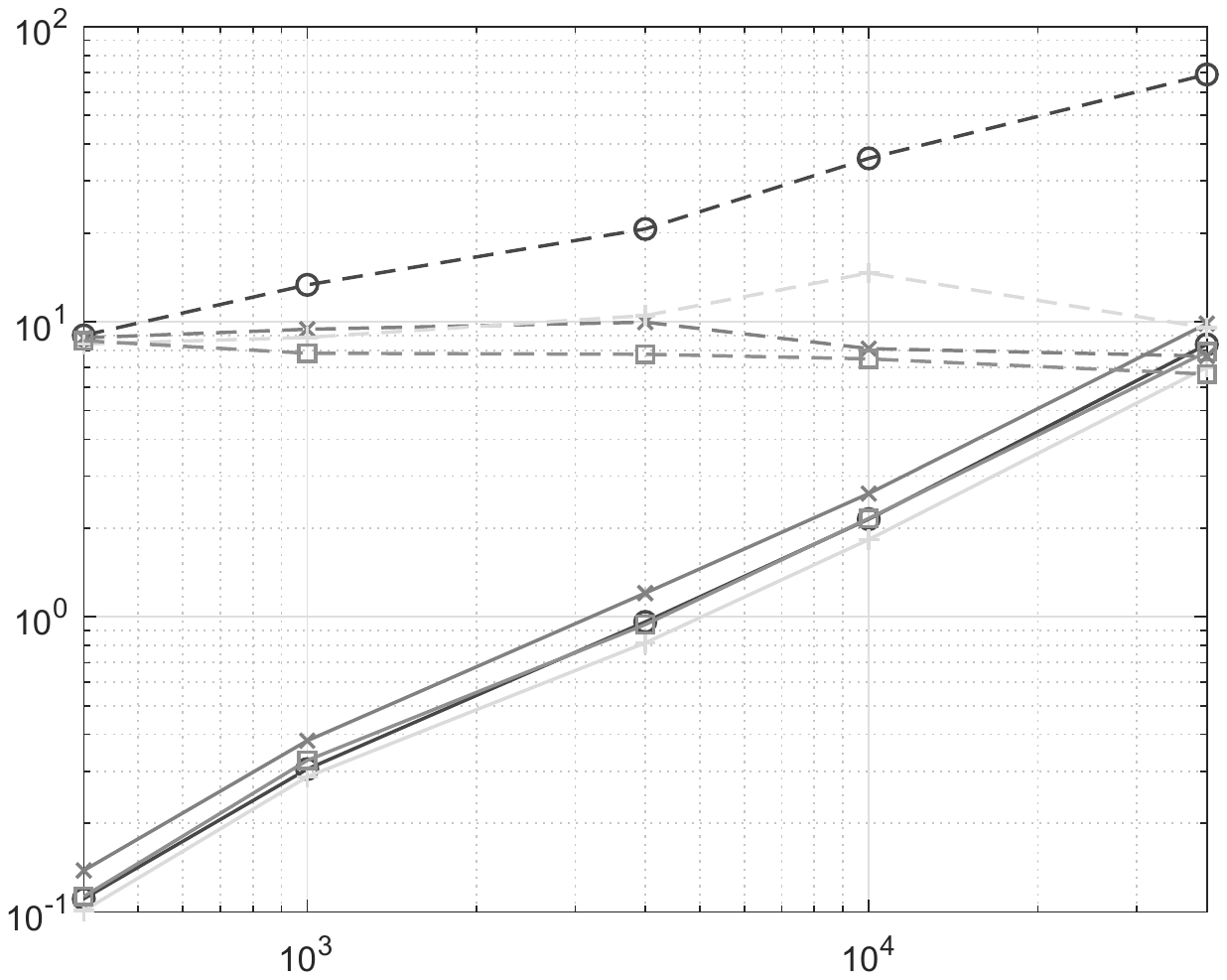}\\
Mean estimation error $\|\widehat{x}_k(\omega^k)-x\|_2$ &CPU time (sec)\end{tabular}
\end{center}\caption{\label{fig:mono} Numerical results: mean errors and CPU times for SA (solid lines) and SAA estimates (dashed lines).
$o$ -- case A (logistic link), x -- case B (linear link), {\small $+$} -- case C (hinge function), {\small $\square$} -- case D (ramp sigmoid).}
\end{figure}
}
Note that the cpu time for SA includes both tuning and execution stages. The conclusion from these experiments is that as far as estimation quality is concerned, the SAA estimate marginally outperforms the SA, while being significantly more time consuming. Note also that the observed in our experiments dependence of recovery errors on $K$ is consistent with the convergence rate $O(1/\sqrt{K})$  established by Proposition \ref{logitprop}.
\section{``Single-observation'' case}\label{logitsingleobs}
Let us look at the special case of the estimation problem where the sequence $\eta_1,...,\eta_K$ of regressors in (\ref{logiteq1})
is deterministic. At the first glance, this situation goes beyond our setup, where the regressors should be i.i.d. drawn from some distribution $Q$. We can, however, circumvent this ``contradiction'' by saying that we are now in the {\sl single-observation case} with the regressor being the matrix $[\eta_1,...,\eta_K]$ and $Q$ being a degenerate distribution supported at a singleton. Specifically, consider the case where our observation is
\begin{equation}\label{logiteq101}
\omega=(\eta,y)\in\bR^{n\times mK}\times \bR^{mK}
\end{equation}
($m,n,K$ are given positive integers), and the distribution $P_x$ of observation stemming from a signal $x\in\bR^n$ is as follows:
\begin{itemize}
\item $\eta$ is a given independent of $x$ deterministic matrix;
\item $y$ is random, and the distribution of $y$ induced by $P_x$ is with mean $\phi(\eta^T x)$, where $\phi:\bR^{mK}\to\bR^{mK}$ is a   given mapping.
\end{itemize}
As an instructive example connecting our current setup with the previous one, consider the case where $\eta=[\eta_1,...,\eta_K]$ with $n\times m$ deterministic ``individual regressors'' $\eta_k$, $y=[y_1;...;y_K]$ with random ``individual labels'' $y_k\in\bR^m$ conditionally independent, given $x$, across $k$, and such that the induced by $x$ expectations of $y_k$ are $f(\eta_k^Tx)$ for some $f:\bR^m\to\bR^m$. We  set $\phi([u_1;...;u_K])=[f(u_1);...;f(u_K)]$.  The resulting ``single observation'' model is a natural analogy of the $K$-observation model considered so far, the only difference being that the individual regressors now form a fixed deterministic sequence rather than being a sample of some random matrix.
\par
Same as {everywhere in this paper,} our goal is to use observation (\ref{logiteq101}) to recover the (unknown) signal $x$ underlying, as explained above, the distribution of the observation. Formally, we are now in the case $K=1$ of our previous recovery problem where $Q$ is supported on a singleton $\{\eta\}$ and can use the constructions developed so far. Specifically,
\begin{itemize}
\item The vector field $F(z)$ associated with our problem (it used to be {\small$\bE_{\eta\sim Q}\{\eta f(\eta^Tz)\}$})  is
$$
F(z)=\eta \phi(\eta^Tz),
$$
and the vector field
$
G(z)=F(z)-F(x),
$
$x$ being the signal underlying observation (\ref{logiteq101}), is
$$
G(z)=\bE_{(\eta,y)\sim P_x}\{F(z)-\eta y\}
$$
(cf. (\ref{logiteq44})). Same as before, the signal to recover is a zero of the latter field. Note that now the vector field $F(z)$ is observable, 
and the vector field $G$ still is the expectation, over $P_x$, of an observable vector field:
$$
G(z)=\bE_{(\eta,y)\sim P_x} \{\underbrace{\eta\phi(\eta^Tz)-\eta y}_{G_y(z)}\},
$$
cf. Lemma \ref{logitobs1}.
\item Assumptions \textbf{A.1-2} now read
\medskip\par\noindent\textbf{A.1$^\prime$} The vector field $\phi(\cdot):\bR^{mK}\to\bR^{mK}$ is continuous and monotone, so that $F(\cdot)$ is continuous and monotone as well,
\medskip\par\noindent\textbf{A.2$^\prime$} $\cX$ is a nonempty compact convex set, and $F$ is strongly monotone, with modulus $\varkappa>0$, on $\cX$.
\end{itemize}
A simple sufficient condition for the validity of the above monotonicity assumptions is positive definiteness of the matrix $\eta\eta^T$ plus strong monotonicity of $\phi$ on every bounded set.
\begin{itemize}
\item
For our present purposes, it is convenient to reformulate assumption \textbf{A.3} in the following equivalent form:
\medskip
\par\noindent
\textbf{A.3$^\prime$} For properly selected $\sigma\geq0$ and every $x\in\cX$ it holds
$$
\bE_{(\eta,y)\sim P_x}\{\|\eta [y-\phi(\eta^Tx)]\|_2^2\}\leq \sigma^2.
$$
\end{itemize}
In the present setting, the SAA estimate $\widehat{x}(y)$ is the unique weak solution to $\VI(G_y,\cX)$, and we can easily quantify the quality of this estimate:
{\begin{proposition}\label{logitpropSAA} In the situation in question, let Assumptions \textbf{A.1$^\prime$-3$^\prime$} hold. Then for every $x\in\cX$ and every realization $(\eta,y)$ of induced by $x$ observation {\rm (\ref{logiteq101})} one has
\begin{equation}\label{logiteq102}
\|\widehat{x}(y)-x\|_2\leq \varkappa^{-1}\|\underbrace{\eta[y-\phi(\eta^Tx)]}_{\Delta(x,y)}\|_2,
\end{equation}whence also
\begin{equation}\label{logiteq104}
\bE_{(\eta,y)\sim P_x}\{\|\widehat{x}(y)-x\|_2^2\}\leq \sigma^2/\varkappa^2.
\end{equation}
\end{proposition}}\noindent
\textbf{Proof.} Let $x\in\cX$ be the signal underlying observation (\ref{logiteq101}), and $G(z)=F(z)-F(x)$ be the associated vector field $G$. We have
$$
G_y(z)=F(z)-\eta y=F(z)-F(x)+[F(x)-\eta y]=G(z)-\eta[y-\phi(\eta^Tx)]=G(z)-\Delta(x,y).
$$
For $y$ fixed, $\bar{z}=\widehat{x}(y)$ is the weak, and therefore the strong (since $G_y(\cdot)$ is continuous) solution to $\VI(G_y,\cX)$, implying, due to $x\in\cX$, that
$$
0\leq G_y^T(\bar{z})[x-\bar{z}]=G^T(\bar{z})[x-\bar{z}]-\Delta^T(x,y)[x-\bar{z}],
$$
whence
$$
-G^T(\bar{z})[x-\bar{z}]\leq-\Delta^T(x,y)[x-\bar{z}].
$$
Besides this, $G(x)=0$, whence $G^T(x)[x-\bar{z}]=0$, and we arrive at
$$
[G(x)-G(\bar{z})]^T[x-\bar{z}]\leq -\Delta^T(x,y)[x-\bar{z}],
$$
whence also
$$
\varkappa\|x-\bar{z}\|_2^2\leq  -\Delta^T(x,y)[x-\bar{z}]
$$
(recall that $G$, along with $F$, is strongly monotone with modulus $\varkappa$ on $\cX$ and $x,\bar{z}\in\cX$). Applying the Cauchy inequality, we arrive at (\ref{logiteq102}). \qed
\par\noindent
\textbf{Example.} Consider the case where $m=1$, $\phi$ is strongly monotone, with modulus $\varkappa_\phi>0$, on the entire $\bR^K$, and $\eta$ in (\ref{logiteq101}) is drawn from a ``Gaussian ensemble'' -- the columns $\eta_k$ of the $n\times K$ matrix $\eta$ are independent $\cN(0,I_n)$-random vectors. Assume also that the observation noise is Gaussian:
$$
y=\phi(\eta^Tx)+\lambda\xi,\;\;\xi\sim\cN(0,I_K).
$$
It is well known that 
 as $K/n\to\infty$, the minimal singular value of the $n\times n$ matrix $\eta\eta^T$ is at least $O(1)K$ with overwhelming probability, implying that when $K/n\gg 1$, the typical modulus of strong monotonicity  of $F(\cdot)$ is $\varkappa\geq O(1)K\varkappa_\phi$. Furthermore, in our situation, as $K/n\to\infty$, the Frobenius norm of $\eta$ with overwhelming probability is at most $O(1)\sqrt{nK}$. In other words, when $K/n$ is large, a ``typical'' recovery problem from the just described ensemble satisfies the premise of Proposition \ref{logitpropSAA} with $\varkappa=O(1)K\varkappa_\phi$ and $\sigma^2=O(\lambda^2nK)$. As a result,
(\ref{logiteq104}) reads
$$
\bE_{(\eta,y)\sim P_x}\{\|\widehat{x}(y)-x\|_2^2\}\leq O(1){\lambda^2n\over\varkappa_\phi^2 K}.\eqno{[K\gg n]}
$$
It is well known that in the standard case of linear regression, where $\phi(x)=\varkappa_\phi x$, the resulting bound is near-optimal, provided $\cX$ is large enough.
\bigskip\par
\begin{figure}[h]
\begin{center}
\begin{tabular}{cc}
\includegraphics[width=0.45\textwidth]{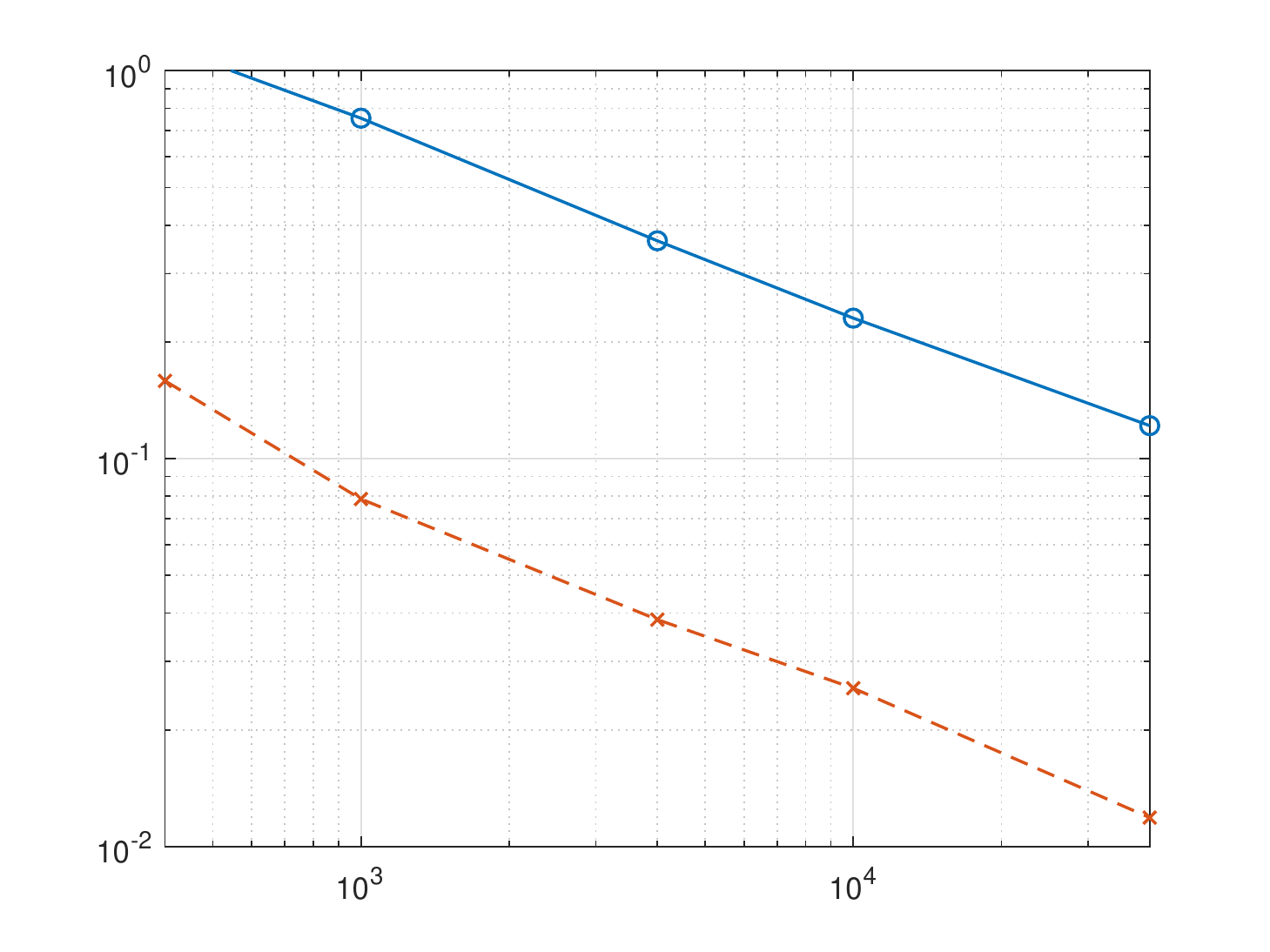}&
\includegraphics[width=0.45\textwidth]{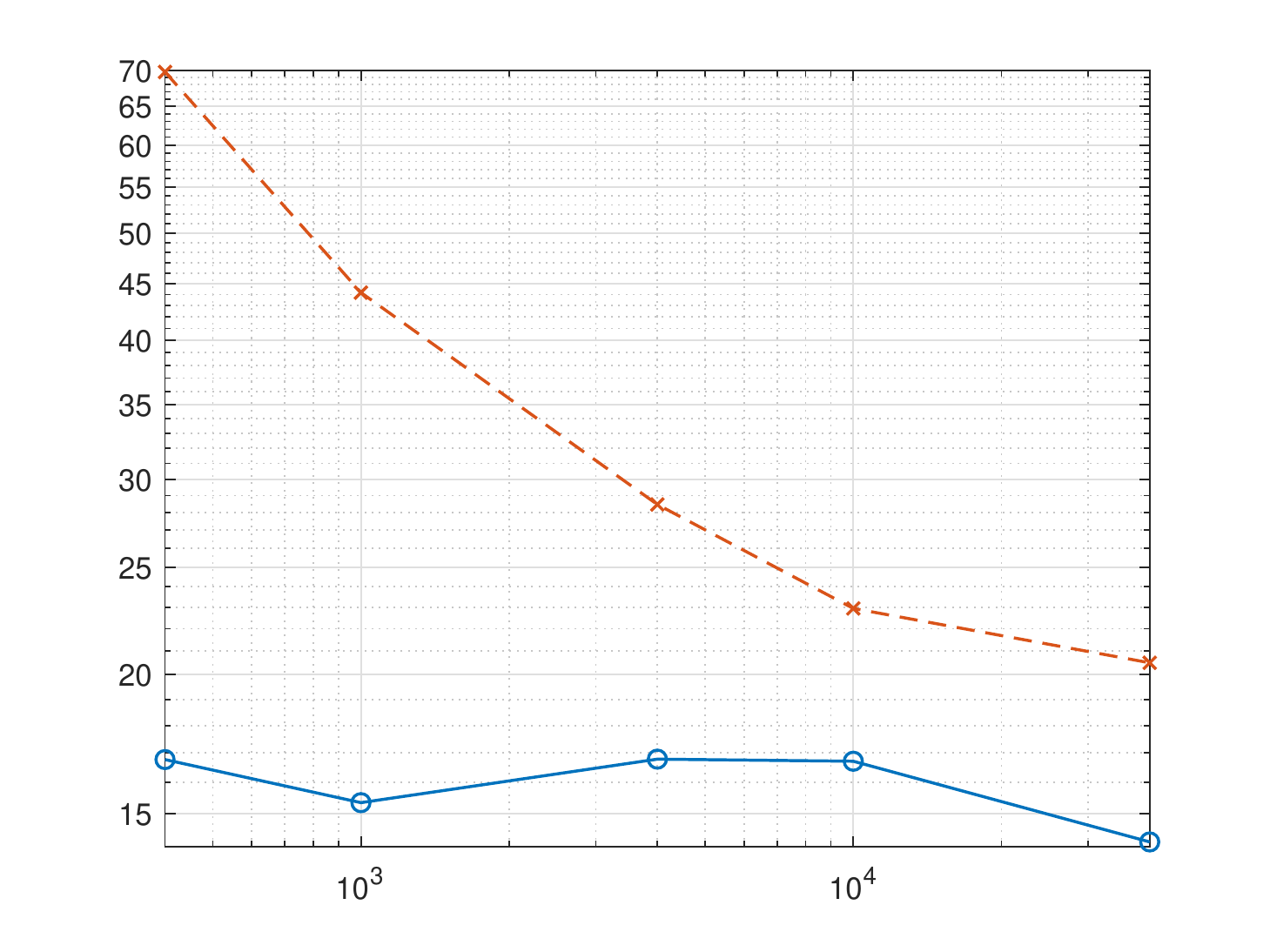}\\
Mean estimation error $\|\widehat{x}_k(\omega^k)-x\|_2$ &CPU time (sec)\end{tabular}
\end{center}\caption{\label{logitfig} Mean errors and CPU times for standard deviation $\lambda=1$ (dashed line) and $\lambda=0.1$ (solid line).}
\end{figure}

\noindent\textbf{Numerical illustration:} in the situation described in the example above, we set $m=1$, $n=100$ and use
$$
\phi(u)=\arctan[u]:=[\arctan(u_1);...;\arctan(u_K)]:\;\bR^K\to\bR^K.
$$
The set $\cX$ is the unit ball $\{x\in\bR^n:\|x\|_2\leq1\}$. In a particular experiment, $\eta$ is chosen at random from the Gaussian ensemble as described above, and signal $x\in\cX$ underlying observation (\ref{logiteq101}) is drawn at
random; the observation noise $y-\phi(\eta^Tx)$ is $\cN(0,\lambda^2 I_K)$. We ran 10 simulations for each combination of the samples size and noise variance $\lambda^2$; the  results are presented 
in Figure \ref{logitfig}.



\begin{thebibliography}{10}

\bibitem{aiserman1964theoretical}
M.~Aiserman, E.~M. Braverman, and L.~Rozonoer.
\newblock Theoretical foundations of the potential function method in pattern
  recognition.
\newblock {\em Avtomat. i Telemeh}, 25:917--936, 1964.

\bibitem{aizerman1970method}
M.~Aizerman, E.~Braverman, and L.~Rozonoer.
\newblock {\em Method of potential functions in the theory of learning
  machines}.
\newblock Nauka, Moscow, 1970.

\bibitem{barndorff1978information}
O.~Barndorff-Nielsen.
\newblock Information and exponential families in statistical theory.
\newblock 1978.

\bibitem{block1962perceptron}
H.-D. Block.
\newblock The perceptron: A model for brain functioning. i.
\newblock {\em Reviews of Modern Physics}, 34(1):123, 1962.

\bibitem{bousquet2004introduction}
O.~Bousquet, S.~Boucheron, and G.~Lugosi.
\newblock Introduction to statistical learning theory.
\newblock In {\em Advanced lectures on machine learning}, pages 169--207.
  Springer, 2004.

\bibitem{bousquet2002stability}
O.~Bousquet and A.~Elisseeff.
\newblock Stability and generalization.
\newblock {\em Journal of machine learning research}, 2(Mar):499--526, 2002.

\bibitem{devroye2013probabilistic}
L.~Devroye, L.~Gy{\"o}rfi, and G.~Lugosi.
\newblock {\em A probabilistic theory of pattern recognition}, volume~31.
\newblock Springer Science \& Business Media, 2013.

\bibitem{devyaterikov1967iterative}
I.~Devyaterikov, A.~Propoi, and Y.~Z. Tsypkin.
\newblock Iterative learning algorithms for pattern recognition.
\newblock {\em Automation and Remote Control}, 28:122--132, 1967.

\bibitem{feigin1981conditional}
P.~D. Feigin et~al.
\newblock Conditional exponential families and a representation theorem for
  asympotic inference.
\newblock {\em The Annals of Statistics}, 9(3):597--603, 1981.

\bibitem{helmbold1995weak}
D.~P. Helmbold and M.~K. Warmuth.
\newblock On weak learning.
\newblock {\em Journal of Computer and System Sciences}, 50(3):551--573, 1995.

\bibitem{KieferSA}
J.~Kiefer and J.~Wolfowitz.
\newblock Stochastic estimation of the maximum of a regression function.
\newblock {\em The Annals of Mathematical Statistics}, 23(3), 1952.

\bibitem{mccullagh1989generalized}
P.~McCullagh and J.~A. Nelder.
\newblock {\em Generalized Linear Models}, volume~37.
\newblock CRC Press, 1989.

\bibitem{nelder1972generalized}
J.~A. Nelder and R.~W. Wedderburn.
\newblock Generalized linear models.
\newblock {\em Journal of the Royal Statistical Society: Series A (General)},
  135(3):370--384, 1972.

\bibitem{NOR}
A.~Nemirovski, S.~Onn, and R.~U.
\newblock Accuracy certificates for computational problems with convex
  structure.
\newblock {\em Mathematics of Operations Research}, 35(1):52--78, 2010.

\bibitem{rakhlin2005stability}
A.~Rakhlin, S.~Mukherjee, and T.~Poggio.
\newblock Stability results in learning theory.
\newblock {\em Analysis and Applications}, 3(04):397--417, 2005.

\bibitem{robbins1951stochastic}
H.~Robbins and S.~Monro.
\newblock A stochastic approximation method.
\newblock {\em The annals of mathematical statistics}, pages 400--407, 1951.

\bibitem{rosenblatt1958perceptron}
F.~Rosenblatt.
\newblock The perceptron: a probabilistic model for information storage and
  organization in the brain.
\newblock {\em Psychological review}, 65(6):386, 1958.

\bibitem{shalev2009stochastic}
S.~Shalev-Shwartz, O.~Shamir, N.~Srebro, and K.~Sridharan.
\newblock Stochastic convex optimization.
\newblock In {\em COLT}, 2009.

\bibitem{ShapiroSAA}
A.~Shapiro, D.~Dentcheva, and Ruszczy\'{n}ski.
\newblock {\em Lectures on Stochastic Programming: Modeling and Theory, Second
  Edition}.
\newblock SIAM, 2014.

\bibitem{sridharan2009fast}
K.~Sridharan, S.~Shalev-Shwartz, and N.~Srebro.
\newblock Fast rates for regularized objectives.
\newblock In {\em Advances in Neural Information Processing Systems}, pages
  1545--1552, 2009.

\bibitem{wolfowitz1952stochastic}
J.~Wolfowitz et~al.
\newblock On the stochastic approximation method of {R}obbins and {M}onro.
\newblock {\em The Annals of Mathematical Statistics}, 23(3):457--461, 1952.

\end{thebibliography}
\appendix
\section{Proof of Proposition \ref{logitprop}}
We start by observing that $z_k$ are deterministic functions of the initial fragments $\omega^{k}=\{\omega_t,1\leq t\leq k\}\sim \underbrace{P_x\times...\times P_x}_{P^k_x}$ of our sequence of observations $\omega^K=\{\omega_k=(\eta_k,y_k),1\leq k\leq K\}$:
$z_k=Z_k(\omega^{k})$.
Let us set
$$
D_k(\omega^{k})=\mbox{$1\over 2$} \|Z_k(\omega^{k})-x\|_2^2=\half \|z_k-x\|_2^2,\quad
d_k=\bE_{\omega^{k}\sim P_x^k}\{D_k(\omega^{k})\},
$$
where $x\in\cX$ is the signal underlying observations (\ref{logiteq1}).
Note that, as it is well known, the metric projection onto a closed convex set $\cX$ {is contracting}:
$$
\forall (z\in\bR^n,u\in\cX):\;\|\Proj_\cX[z]-u\|_2\leq \|z-u\|_2.$$
 Consequently, for $1\leq k
\leq K$ it holds
\bse
D_{k}(\omega^k)&=&\half\|\Proj_{\cX}[z_{k-1}-\gamma_k G_{\omega_k}(z_{k-1})-x]\|_2^2\\
&\leq& \half\|z_{k-1}-\gamma_k G_{\omega_k}(z_{k-1})-x\|_2^2\\
&=&\half\|z_{k-1}-x\|_2^2-\gamma_kG^T_{\omega_k}(z_{k-1})(z_{k-1}-x)+\half \gamma_k^2\|G_{\omega_k}(z_{k-1})\|_2^2.
\ese
Taking expectations w.r.t. $\omega^k\sim {P_x^k}$ of both sides of the resulting inequality and keeping in mind relations
(\ref{logiteq5}) along with the fact that $z_{k-1}\in \cX$, we get
\be
d_{k}\leq d_{k-1}-\gamma_k\bE_{\omega^{k-1}\sim P_x^{k-1}}\left\{G(z_{k-1})^T(z_{k-1}-x)\right\}+2\gamma_k^2M^2.
\ee{logiteq12}
Recalling that we are in the case where $G$ is strongly monotone on $\cX$ with modulus $\varkappa>0$,  $x$ is the weak solution $\VI(G,\cX)$, and $z_{k-1}$ takes values in $\cX$, invoking (\ref{logiteq8}), the expectation in (\ref{logiteq12}) is at least $2\varkappa d_k$, and we arrive at the
relation
\begin{equation}\label{logit24}
d_{k}\leq (1-2\varkappa\gamma_k)d_{k-1}+2\gamma_k^2M^2.
\end{equation}
We put
$$
S={2M^2\over\varkappa^2},\quad \gamma_k={\varkappa S\over 4M^2(k+1)}={1\over\varkappa (k+1)};
$$
note that $\gamma_k$ are exactly the stepsizes (\ref{logiteq50}).
Let us verify by induction in $k$ that for $k=0,1,...,K$ it holds
$$
d_k\leq {(k+1)^{-1}S}.\eqno{(*_k)}
$$
\textbf{Base $k=0$.} Let $D$ stand for the $\|\cdot\|_2$-diameter of $\cX$, and $z_\pm\in\cZ$ be such that $\|z_+-z_-\|_2=D$. By  (\ref{logiteq5}) we have $\|F(z)\|_2\leq M$ for all $z\in\cX$, and by strong monotonicity of $G(\cdot)$ on $\cX$ we have
$$
[G(z_+)-G(z_-)]^T[z_+-z_-]=[F(z_+)-F(z_-)][z_+-z_-]\geq \varkappa \|z_+-z_-\|_2^2=\varkappa D^2;
$$
By Cauchy inequality, the left hand side in the concluding $\geq$ is at most $2MD$, and we get
$$
D\leq {2M\over\varkappa},
$$
whence $S\geq D^2/2$. On the other hand, due to the origin of $d_0$ we have $d_0\leq D^2/2$. Thus, $(*_0)$ holds true.\\
\textbf{Inductive step $(*_{k-1})\Rightarrow (*_{k})$.} Now assume that $(*_{k-1})$ holds true for some $k$, $1\leq k\leq K$, and let us prove that $(*_{k})$ holds true as well. Observe that $\varkappa\gamma_k=(k+1)^{-1}\leq 1/2$, so that
\bse
d_{k}&\leq& d_{k-1}(1-2\varkappa\gamma_k)+2\gamma_k^2M^2\hbox{\ [by (\ref{logit24})]}\\
&\leq& {S\over k}(1-2\varkappa \gamma_k)+ 2\gamma_k^2M^2\hbox{\ [by $(*_{k-1})$ and due to $\varkappa\gamma_k\leq1/2$]}\\
&=&{S\over k}\left(1-{2\over k+1}\right)+{S\over (k+1)^2}={S\over k+1}\left({k-1\over k}+{1\over k+1}\right)\leq  {S\over k+1},
\ese
so that $(*_{k})$ hods true. Induction is complete. It remains to note that by definition of $d_k$ we have $d_k={1\over 2}\bE\{\|\widehat{x}_k-x\|_2^2\}$. \qed
\end{document}